\documentclass{amsart}
\usepackage{amscd}
\usepackage{amsmath}

\newtheorem{theorem}{Theorem}[section]
\newtheorem{lemma}[theorem]{Lemma}

\theoremstyle{definition}

\theoremstyle{remark}

\numberwithin{equation}{section}

\begin{document}

\def\reals{\text{{\rm I \kern -4.7pt R}}}
\def\rmap{\mbox{\bf r}}
\def\lieg{{\mathfrak g}}
\def\liek{{\mathfrak  k}}
\def\q{\mathfrak{\bf q}}
\def\barM{
\hbox{\kern 2.3 true pt \vbox{\hrule width 8.5  true pt height .3
true pt \kern .9 true pt \hbox{\kern -2.3 true pt $M$}}}}
\def\Vert{\hbox{\it Vert} }
\def\Vertbundle{\hbox{{\it Vert}}\, M}
\def\Vertbundlept#1{\Vert_{#1}\kern-.5pt M }
\def\bGamma{{ \boldsymbol\Gamma }}
\def\TM{T\kern -.5pt M}
\def\mychi{{\hbox{\raise 2 pt \hbox{$\chi$}}}}
\def\Ad{\text{\rm Ad}}
\def\LD{\mathcal{ L}}

\def\barM{
\hbox{\kern 2.3 true pt \vbox{\hrule width 8.5  true pt height .3
true pt \kern .9 true pt \hbox{\kern -2.3 true pt $M$}}}}
\def\balpha{{ \bar\alpha }}
\def\barU{
\hbox{ \vbox{\hrule width 7 true pt height .3 true pt \kern .9
true pt \hbox{\kern-1 pt $U$}}}}

\title{ A Co-chain map for the G-invariant
\\
de Rham complex.}

\author{ I.M. Anderson, M.E. Fels }
\address{ Department of Mathematics and Statistics, \\
Utah State University, Logan Utah,USA, 84322 \\ E-mail:
anderson@math.usu.edu, fels@math.usu.edu }

\maketitle


\section{Introduction}

In this note we characterize the Lie group actions for which there
exists, at least locally, an evaluation map that defines a cochain
map from the differential complex of invariant forms on a manifold
to the De Rham complex for the quotient. This problem is motivated
by the principle of symmetric criticality \cite{and3}.

Before giving any specific definitions we would like to illustrate
the notion of such an evaluation map with a simple example.
Consider the two dimensional Abelian Lie group $G=\reals^2$ with
coordinates $(a,b)$ acting  on $\reals^3$ by
$$
(a,b)* (x,y,z) = (x,y+a,z+b) .
$$
If $ \alpha \in \Omega^2(\reals^3)^G$ and $\nu \in
\Omega^3(\reals^3)^G$, where we use the convention that a group
superscript denotes the invariants of the group, then $\alpha$ and
$\nu$ are necessarily of the form
$$
\alpha = a(x) dx\wedge dy  + b(x) dx \wedge dz + c(x) dy \wedge dz
\ \ {\rm and} \quad \nu = A(x) dx\wedge dy \wedge dz.
$$
The Lie algebra of infinitesimal generators of this action of $G$
is generated by $\{\partial_y,\partial_z\}$ and it is easy check
that evaluation on the generators
$$
\alpha(\partial_y,\partial_z ) = c(x) \ , \quad
 \nu(\partial_y,\partial_z,-)
= A(x) dx
$$
defines a cochain map from $\Omega^*(\reals^3)^G$ to
$\Omega^{*-2}(\reals)$, that is,
$$
(d \alpha)(\partial_y,\partial_z,-) =
d(\alpha(\partial_y,\partial_z)) = c(x)'  dx.
$$

As we shall see, not all group actions admit cochain evaluation
maps.

\section{Lie group actions and
invariant vector fields }

Let $G$ be a $p$-dimensional Lie group which acts effectively on
an $n$-dimensional manifold $M$ with multiplication map
$\mu:G\times M \to M$. We write $gx$ instead of $\mu(g,x)$. For
$x\in M$ and $g\in G$, we define $ \mu_x :G\to M $ and $ \mu_g:M \to M$
to be the maps
$$
\mu_x(g) = \mu_g(x)=g x.
$$
For any $g\in G$, $\mu_g$ is a   diffeomorphism of $M$. We let $G_x$  denote the
isotropy subgroup of $G $ at $x$,
$$
    G_x=\{\,  g\in G \, | \, g x = x \,\} .
$$
For each $x\in M$, the map $\tilde \mu_x : G/G_x \to M$ given by
$\tilde\mu_x ([g]) = g x $ is a one-to-one immersion which is also
$G$ equivariant with respect to the canonical action of $G$ on the
coset space $G/G_x$.

The Lie algebra $\lieg$ of the Lie group $G$ is the Lie algebra of
{\it right } invariant  vector fields on $G$. The action of $G$ on
$M$ induces a Lie algebra homomorphism  $ \rmap : \lieg \to
\mathcal{ X}(M)$  of $\lieg $ to the vector fields on $M$ whose
image is the Lie algebra of the infinitesimal generators of the
action of $G$ on $M$ \cite{olv1}. We write $\Gamma =
\rmap(\lieg)$. Because the action of $G$ on $M$ is assumed
effective, the map $\rmap$ is injective.
 Let $ \bGamma \subset TM $ denote the (integrable)
distribution generated by $ \Gamma $.

The action of the Lie group $G$ on $M$ is said to be {\it regular}
if the space of orbits is a manifold $\barM =M/G$ such that the
quotient map
$$
    \q : M \to \barM
$$
is a submersion. We will assume from here on that all actions are
regular. For regular actions the orbits all have the same
dimension which we assume to be $q$ and so the isotropy subgroup
$G_x$, for any $x\in M$, will have dimension $p-q$. Let
$\Vertbundle \to M$ be the sub-bundle of $\q$ vertical vectors in
$\TM$, so $\Vertbundle = \ker \q_* = \bGamma$. We also have the
important property $(\tilde \mu_x)_*(T_{[e]}G/G_x)=
\Vertbundlept{x}$.

The action of $G$ on $M$ defines an action of $G$ on $TM$
using the differential
\begin{equation}
(\mu_g)_* : TM \to TM .
\label{eq:dg}
\end{equation}
For each $g \in G_x$, equation (\ref{eq:dg}) gives
$$
    (\mu_g)_* : T_{\kern- .8pt x} \kern -.5pt M \to  T_{\kern -.8ptx} \kern -.5pt  M
$$
which defines the linear isotropy representation of $G_x$ on the
tangent space $T_{\kern -.8pt x} \kern -.5pt M$.

Suppose now that $X$ is a $G$ invariant vector field, that is,
\begin{equation}
(\mu_g)_* X_x = X_{gx} .
\label{eq:isor}
\end{equation}
If $g \in G_x$, then equation (\ref{eq:isor}) implies that
\begin{equation}
X_x \in (T_xM)^{G_x}.
\label{eq:isfp}
\end{equation}
This observation leads us to define the following subset of $TM$,
$$
\kappa(TM) = \bigcup_{x\in M}\kappa(T_xM)\, , \ {\rm where} \quad
\kappa(T_xM) = (T_xM)^{G_x}.
$$
Equation (\ref{eq:isfp}) implies that every $G$ invariant vector
field $X$ takes values in the subset $\kappa(TM) \subset TM$.

Since $\q\circ \mu_g =\q$, the action of $G$ on $TM$ restricts to
an action on $\Vertbundle$ and the linear isotropy representation
of $G_x$ also restricts to a representation on vertical vectors
$$
    (\mu_g)_* : \Vertbundlept{x} \to \Vertbundlept{x}, \quad g \in
    G_x.
$$
Thus a $G$ invariant vertical vector field takes values in the set
$$
\kappa(\Vertbundle) = \bigcup_{x\in M}\kappa(\Vertbundlept{x})\, ,
\ {\rm where} \quad \kappa(\Vertbundlept{x}) =
(\Vertbundlept{x})^{G_x}.
$$

In the next theorem we give  conditions which guarantee the
existence of invariant vector fields. This is a special case of
the general construction given in \cite{and1} or on
 p. 657 in \cite{and2}.

\begin{theorem} If $\kappa(TM) \subset TM$ is a vector sub-bundle, then for
each $x\in M$ and $Y \in \kappa(T_xM)$ there exists a $G$
invariant vector field $X$ on  $M$ such that $X_x=Y$. The
analogous statement holds for $G$ invariant vertical vector fields
if $\kappa(\Vertbundle)\subset
 \Vertbundle$ is a vector sub-bundle.
\end{theorem}

\noindent
{\bf Remark 2.1} For the rest of this article we assume that all group
actions are regular and that $\kappa(TM)$, and $\kappa(\Vertbundle)$
are bundles.

\section{ Lie algebra cohomology}

Given a Lie group $G$ and a Lie subgroup $K\subset G$,  with
corresponding Lie algebras $\liek\subset \lieg$,  define the
vector space of $K$   relative forms on $\lieg $ by
$$
    A ^r (\lieg ,  K ) = \bigr\{\, \alpha \in A ^r (\lieg) \,| \, \iota_v \alpha = 0
    , \ \forall \ v \in \liek \ \text{and}\  \Ad^* g \cdot \alpha
= \alpha
                \ , \ \forall  \ g \in K \, \bigl\},
$$
where $A^r(\lieg)$ are the alternating $r$-forms on $\lieg$ and
 $ \Ad^* $ denotes the co-adjoint representation of $G$ on
$A^r(\lieg)$.

The usual   exterior derivative $d$ on $ A^*(\lieg)$ restricts to
make $A^*(\lieg, K)$  a differential complex whose cohomology is
denoted by
 $H^*(\lieg, K)$, the Lie algebra cohomology of $\lieg$ relative to the
subgroup $K$.

If $K\subset G$ is a closed Lie subgroup, let $H^*(\Omega^*(G/K)^G)$ be the
$d$-cohomology of the $G$ invariant forms on $G/K$.

\begin{lemma} If $K \subset G$ is closed, then
$\Omega^r(G/K)^G \simeq  A ^r (\lieg ,  K )$ and $H^r(
\Omega^*(G/K)^G) \simeq H^r(\lieg, K)$.
\end{lemma}
See Theorem 13.1 in \cite{chev} for a proof of this Lemma. It is
well-known \cite{spiv}, that if $G$ is connected and compact and
$K$ closed, then $H^*(\lieg, K) $ computes the De Rham cohomology
of the homogeneous space $G/K$.

It is useful to note that if $K_2=gK_1g^{-1}$ are conjugate
subgroups of $G$ then $\Ad (g)$ induces an isomorphism
\begin{equation}
    A^*(\lieg, K_1) = A^*(\lieg, K_2).
\label{eq:conj}
\end{equation}

\smallskip
\noindent {\bf Example:}  Consider the two sphere $S^2$ and the
projective plane $\reals {\rm P}^2$ as the homogeneous spaces
$SO(3)/SO(2)$ and $SO(3)/O(2)$. Letting $ X_1,X_2,X_3$ be a basis
for $so(3)$ with $X_3$ the basis for $so(2)$ (which particular
$so(2)$ is actually irrelevant because of (\ref{eq:conj})) and
letting $\alpha^1,\alpha^2,\alpha^3$ be the dual basis, we find
$$
    A^1( so(3),SO(2)) =  \{ 0 \}
            \quad\text{and}\quad
           A^2(so(3),SO(2)) = \{ \alpha ^1 \wedge \alpha^2 \} .
$$
Therefore $H^2(so(3), SO(2))$ is generated by $\alpha_1\wedge
\alpha_2$. On the other hand, there is  a reflection in $O(2)$
which maps $X_1$ to $-X_1$ and $X_2$  to $X_2$ so that
$$
    A^1( so(3),O(2)) =  \{ 0 \} \quad\text{and}\quad  A^2(so(3),O(2)) = \{0\}
$$
and  therefore $ H^2(so(3), O(2))=0$. Of course, these
computations reflect  the fact that $S^2$ is orientable  whereas
$\reals  {\rm P}^2$ is not.

\section{ A map on the G invariant De Rham complex }

In this section we generalize the evaluation map from the
introduction by studying the problem of defining a map
$$
\rho^k_\mychi : \Omega^{k}(M)^G \to \Omega^{k-q}(\barM)
$$
which shifts  form degree by the orbit dimension $q$ of $G$ on
$M$. To begin, we define $\Lambda_q(\Vertbundle) \to M$ to be the
vector bundle of vertical  $q$-chains on $\Vertbundle$
(alternatively, the bundle of vertical multi-vectors of degree
$q$). Given that the orbits of $G$ have dimension $q$ it follows
that about each point $x\in M$ there exists an open set $U$ and
vector fields $X_1$, $X_2$, \dots, $X_q$ in $\Gamma$ which define
a local frame for $\bGamma|_U =\Vert\, U$. Consequently if
$\mychi$ is a section of $\Lambda_q(\Vertbundle)$ then $\mychi|_U$
can written as
$$
    \mychi|_U= J \, X_1 \wedge X_2  \wedge \cdots \wedge X_q ,
$$
where $J \in C^\infty(U)$. The  action of $G$ on $\Vertbundle$
described in section 2, induces an action of $G$ on
$\Lambda_q(\Vertbundle) $.

Given  a $G$ invariant  $q$-chain $\mychi :M \to
\Lambda_q(\Vertbundle)$, we now define a map $
 \iota_\mychi : \Omega^k(M)  \to \Omega^{k-q}_{\rm sb}(M)$
where
$$
\Omega^*_{sb}(M) = \{ \omega \in \Omega^*(M) \ | \ \iota_X \omega
= 0 \quad {\rm for \ all} \ X \in \bGamma \}
$$
are the $\q$ semi-basic forms on $M$. The map $\iota_\mychi$ is
defined by setting
$$
(\iota _\mychi \omega)_x(Y_1, Y_2,\dots, Y_{k-q})          =
\omega_x(\mychi_x, Y_1, Y_2,\dots, Y_{k-q})
$$
for  $\omega \in \Omega^k(M)$ and $Y_i \in T_{x}M$. If $ \omega
\in \Omega^k(M)^G$ then $\iota_\mychi \omega$ is $\q$ semi-basic,
and since $\mychi$ is $G$ invariant, $\iota_\mychi \omega$ is $G$
invariant and so $G$ basic.  By this last statement $\iota_\mychi
\omega \in \Omega^{k-q}_{sb}(M)^G$, and therefore by Lemma A.3 in
\cite{and0}, we find there exists a unique $(k-q)$-form
$\overline{\iota_\mychi \omega}$ on $\barM$ satisfying
$\q^*(\overline{\iota_\mychi \omega}) = \iota_\mychi \omega$. The
sought after evaluation map $\rho_\mychi$ is then defined by
\begin{equation}
           \rho^k_\mychi(\omega) =  (-1)^{(n-k)q}\, \overline{ \iota _\mychi  \omega}.
\label{eq:rhom}
\end{equation}
Note that for each invariant $\mychi$ we have a map $\rho_\mychi$.

\begin{theorem} If there exists a non-vanishing $G$
invariant  vertical  $q$-chain $\mychi$ on $M$, then
\begin{equation}
    A^q(\lieg , G_x) \neq 0\  \quad { for \ all} \ x\in M.
\end{equation}
Conversely, if for each $x\in M$, $A^q(\lieg , G_x) \neq 0$ then
about each $x_0 \in M$ there exists a $G$ invariant open set $U$
and non-vanishing $G$ invariant vertical $q$-chain $\mychi$ on
$U$.
\end{theorem}

\begin{proof} Let $\mychi$ be a non-vanishing
$G$ invariant vertical $q$-chain. Let $x\in M$ and let $\tilde
\mychi$ be the restriction of $\mychi$ to $G/G_x$, so that
$(\tilde \mu_x )_* \tilde \mychi = \mychi $. By the equivariance
property of $\tilde \mu_x$ the $q$-chain $\tilde \chi $ is $G$
invariant. Now let $\alpha \in \Omega^q(M)$ satisfy
$\alpha(\mychi) = 1 $. The form $\alpha $ is not unique, and it is
not necessarily invariant. We claim the form ${\tilde \mu_x}^*
\alpha $ defines a non-zero element of $\Omega^q(G/G_x)^G $. We
compute
$$
(g^* \tilde \mu_x^* \alpha)(\tilde \mychi) = \alpha ( (\tilde
\mu_x)_* g_* \tilde \mychi ) = \alpha( (\tilde \mu_x)_* \tilde
\mychi) = \alpha (\mychi) = 1 .
$$
Thus $\tilde \mu_x^* \alpha$ is a non-vanishing $G$ invariant form
of top degree on $G/G_x$ and so, by Lemma 3.1, $A^q(\lieg , G_x)
\neq 0$.

\smallskip

We now prove the converse part of the theorem. Let
$$
\kappa(\Lambda_q (Vert M)) = \bigcup_{x\in M} \kappa( \Lambda_q
(Vert _xM)),
$$
where $\kappa( \Lambda_q (Vert _xM) )  = (\Lambda_q (Vert
_xM))^{G_x}.$ We shall show that $A^q(\lieg , G_x) \neq 0$ implies
$\kappa(\Lambda_q (Vert M))$ is a line bundle. Then  the existence
of a $G$ invariant $q$-chain is guaranteed (in a similar manner to
Theorem 2.1) by Theorem 1.2 in \cite{and1}.

If $A^q(\lieg , G_x) \neq 0$ then by Lemma 3.1 there exists a
non-vanishing $\tilde \alpha \in \Omega^q(G/G_x)^G$. Let $\tilde
\mychi$ be the invariant $q$-chain defined by $\tilde \alpha(
\tilde \mychi ) = 1$. Then $\chi_x = (\mu_x)_* \tilde \chi _{[e]}
\in \Lambda_q (Vert _xM))^{G_x}$ by the equivariance of $\tilde
\mu_x$, and is non-zero. Thus $\Lambda_q (Vert _xM))^{G_x} =
\Lambda_q (Vert _xM)$ and so $\kappa(\Lambda_q (Vert M)) =
\Lambda_q (Vert M)$ is a line bundle.
\end{proof}

\section{  The cochain condition }

In this section we find necessary and sufficient conditions on the
action of $G$ on $M$ that determine whether we can choose a $G$
invariant $q$-chain $\mychi$  so that the map
$\rho_\mychi:\Omega^*(M)^G \to \Omega^{*-q}(\barM) $ defined in
(\ref{eq:rhom}) is a cochain map, that is,
\begin{equation}
\rho_\mychi(d\,\omega) = d\, \rho_\mychi(\omega). \label{eq:coch}
\end{equation}
Granted that the action of $G$ on $M$ satisfies the conditions in
Remark 2.1, the solution to this problem is given by the following
theorem.

\begin{theorem} If there exists a non-vanishing invariant $q$-chain $\mychi$ such that the map
$\rho_\mychi$ in (\ref{eq:rhom}) defines a cochain map, then $
H^q(\lieg, G_x) \neq 0 $ for all $x \in M.$ Conversely, if
$H^q(\lieg, G_x) \neq 0$ for all $x \in M$ then
 about each $x_0 \in M$ there exists a $G$ invariant open set $U$ and a non-vanishing $G$ invariant
vertical $q$-chain $\mychi$ on $U$ such that
$\rho_\mychi:\Omega^*(U)^G \to \Omega^{*-q}(U/G)$ is a cochain
map.
\end{theorem}

In order to prove this theorem, we need a number of preliminary
results. The first of these is the important observation that the
cochain condition (\ref{eq:coch}), which is a condition that
involves the quotient manifold $\barM$, can be expressed as a
condition entirely on $M$.

\begin{lemma}{\label{lm51}} A $G$ invariant,  vertical $q$-chain
$\mychi$ defines a cochain map $\rho_\mychi$ if and only if
\begin{equation}
\iota_\mychi  d \omega = (-1)^q d ( \iota_\mychi \omega) \qquad
\text{for all} \ \omega \in \Omega^*(M)^G. \label{eq:coup}
\end{equation}
\end{lemma}
\begin{proof} If $\eta$  is any $G$ basic form, then  $d\eta$ is
also $G$ basic. Let $\bar\eta$ be the unique form on $\barM$ such
that $\q^*(\bar \eta) =\eta$. Then, since
$$
\q^*( d \bar\eta) = d \q^*(\bar\eta) = d\eta
$$
the two forms $d\bar\eta$ and $\overline{d\eta}$ pullback by $\q$
to the same form and must therefore  be equal. Since $\mychi$ and
$\omega$ are both $G$ invariant, we can apply this observation to
the $G$ basic form $\iota_\mychi \omega$ to deduce that
$$
    d(\overline{\iota_\mychi \omega} )  = \overline{d(\iota_\mychi \omega)}.
$$
The cochain condition (\ref{eq:coch}) can therefore be expressed
as
\begin{equation}
    (-1)^q\overline{\iota_ \mychi d (\omega)}=  \overline{d(\iota_\mychi\omega)}.
\label{eq:pco}
\end{equation}
But  two $G$ basic forms on $M$  are equal if and only if the
corresponding forms on $\barM$ are equal and so (\ref{eq:pco})
proves the equivalence of (\ref{eq:coch}) with (\ref{eq:coup}).
\end{proof}

\begin{lemma}{\label{lm52}}  If (\ref{eq:coup}) holds for all $G$ invariant
$(n-1)$-forms, then  (\ref{eq:coup})  holds  for all  $G$
invariant $r$-forms, $r \geq q$.
\end{lemma}
\begin{proof} Suppose (\ref{eq:coup}) holds true for all $G$ invariant
$(n-1)$-forms. Let $\omega$ be a $G$ invariant $r$-form, where $ q
\leq r < n-1$. Then, if $ \alpha $  is any $G$ basic
$(n-r-1)$-form, $\omega \wedge \alpha $ is a $G$ invariant
$(n-1)$-form and therefore we can use (\ref{eq:coup}) to write
$$
    \iota_\mychi d (\omega \wedge \alpha) = (-1)^q
            d \left(\iota_\mychi  (\omega \wedge \alpha)  \right) .
$$
Because $\alpha $ (and hence $d \alpha$) is $G$ basic, the
expansion of  both sides of this  equation   gives
$$
(\iota_\mychi  d\omega )\wedge \alpha = (-1)^q  d(\iota_ \mychi
\omega) \wedge \alpha .
$$
 Since $\alpha$ is an arbitrary $G$ basic form and $\iota_\mychi d\omega$ and
$ d( \iota_\mychi \omega)$ are both $G$ basic this implies
$$
    \iota_\mychi d\omega  = (-1)^q d( \iota_\mychi  \omega).
$$
\end{proof}

\begin{lemma}{\label{lm53} } Let  $\mu$  be  a $G$ basic
$(n-q)$-form on a $G$ invariant open set $U$.  Let   $\mychi$ be a
non-vanishing, $G$ invariant, vertical $q$-chain on $U$ and let
$\alpha$ be any $q$-form such that $\alpha(\mychi) = 1$. Then
$$
    \nu = \alpha\wedge \mu
$$
is a $G$ invariant $n$-form on $U$.
\end{lemma}

\begin{proof}  For any $g\in G$, we compute
$$
    [\mu_g^*(\alpha)](\mychi) = \alpha( (\mu_g)_*(\mychi)) =
    \alpha(\mychi)= 1
$$
and therefore
$$
    \iota_\mychi[\mu_g^*(\nu)]= \iota_\mychi[\mu^*_g(\alpha )\wedge \mu] =
             \mu = \iota_\mychi \nu .
$$
This suffices to prove that $\mu_g^*(\nu) = \nu$.
\end{proof}

\begin{lemma}{\label{lm54}} If $\mychi$ is non-vanishing vertical
$q$-chain and $R$ is a $G$ invariant vector field then
$$
    \LD_R\mychi = \lambda_R \mychi,
$$
where $\lambda_R$ is a $G$ invariant function.
\end{lemma}

\begin{proof} Let $ X_1,\dots ,X_q$ be vector fields in $\Gamma$
which form a local basis  for $\Vertbundle$  in some neighborhood
about the point $x$.
 Then $ \mychi = J\, X_1\wedge X_2\wedge\dots \wedge X_q $ and, since
$[\,R, X_i\,] = 0$,
$$
    \LD_R \mychi = R(J)\, X_1\wedge X_2 \wedge \dots\wedge X_q  = \frac{R(J)}{J} \,
 \mychi.
$$
\end{proof}

\begin{theorem} If $\mychi$ is a non-vanishing, $
G$ invariant, vertical $q$-chain, then the map $\rho_\mychi
:\Omega^{*}(M)^G \to \Omega^{*-q}(\barM)$ is a cochain map if and
only  if
\begin{equation}
    \LD_{R} \mychi = 0
\label{eq:ldchi}
\end{equation}
for all $G$ invariant vector fields $R$ on $M$.
\end{theorem}

\begin{proof} We start by assuming (\ref{eq:ldchi}). Then by  Lemma~\ref{lm52}
it suffices to prove (\ref{eq:coup}) for $G$  invariant
$(n-1)$-forms.

Given the non-vanishing $G$ invariant vertical $q$-chain $\mychi$,
let $x\in M$ and use Lemma~\ref{lm53} to construct a non-vanishing
$G$ invariant $n$-form $\nu=\alpha \wedge \mu$ on an invariant
open set $U$ about $x$. Let $\omega \in \Omega^{n-1}(M)^G$, then
restricted to $U$ there exists a unique $G$ invariant vector field
$S$ on $U$ such that
\begin{equation}
    \omega_U =\iota_ S  \nu .
\label{eq:omr}
\end{equation}
Let $R$ be a $G$ invariant vector field on $M$ which agrees with
$S$ on an invariant open set $V\subset U$ of $x$ so that
\begin{equation}
    \omega_V =\iota_ R  \nu.
\label{eq:oms}
\end{equation}
With $\omega_V$  given by (\ref{eq:oms}), we  compute on $V$
\begin{eqnarray*}
   \left[ d(\iota_\mychi \omega)\right]_V & = &d(\iota_\mychi \iota_ R \nu) = (-1)^q d(\iota_R\iota_\mychi\nu)
         = (-1)^q d(\iota_R\mu) \ \ {\rm and}
\\
\left[\iota_\mychi  d\omega \right]_V & = &\iota_\mychi
d(\iota_R\nu) = \iota_\mychi \LD_R(\nu)
            = \LD_R( \mu) - \iota_{\LD_R(\mychi)}  \nu.
\end{eqnarray*}
But it is easy to check that if $\mu$ is a $G$ basic $(n-q)$-form,
then $d\mu = 0$ and  therefore
\begin{equation}
       \left[  d(\iota _\mychi \omega)
           - (-1)^q \iota_\mychi  d\omega \right]_V  =(-1)^q \iota_{\LD_R(\mychi)}   \nu .
\label{eq:ch2}
\end{equation}
Evaluating (\ref{eq:ch2}) at $x\in V$ shows that if
(\ref{eq:ldchi}) holds at $x$ then (\ref{eq:coup}) holds at $x$
for all $G$ invariant $(n-1)$-forms $\omega$. Since our original
point $x\in M$ was arbitrary, equation (\ref{eq:coup}) holds on
$M$.

To prove that (\ref{eq:coup}) implies (\ref{eq:ldchi}) we reverse
the argument above. Let $R$ be a $G$ invariant vector field on $M$
and let $x \in M$. Choose a $G$ basic $(n-q)$-form $\mu$ which
doesn't vanish at $x$. Then the form $\omega = \iota_R \nu$, where
$\nu = \alpha \wedge \mu $ with $\alpha(\mychi)=1$, is a $G$
invariant $n-1$ form on $M$. Equation (\ref{eq:ch2}) (evaluated at
$x$) coupled with Lemma~\ref{lm54} shows that (\ref{eq:coup})
implies (\ref{eq:ldchi}) at $x$. But $x$ was arbitrary so
$(\ref{eq:ldchi})$ holds on $M$.
\end{proof}

We are now in a position to prove Theorem 5.1.

\begin{proof}   We begin the proof by first noting that the condtion $H^q(\lieg,G_x)\neq 0$
is equivalent to the following:

i] For each $x\in M$ there exists a non-vanishing $\tilde \alpha
\in \Omega^q(G/G_x)^G$; and

ii] for all $\tilde \eta \in \Omega^{q-1}(G/G_x)^G$, $d\tilde \eta = 0$.

\smallskip

Suppose there exists a non-vanishing $G$ invariant $q$-chain
$\mychi$ such that $\rho_\mychi$ is a cochain map. Let $\alpha \in
\Omega^q(M)$ with $\alpha(\mychi) = 1$. Then as was shown in
Theorem 4.1, given any $x\in M$,  $\mu_x^* \alpha \in
\Omega^q(G/G_x)^G$ and is non-vanishing.  Thus condition i] is
true.

Let $\tilde \eta \in \Omega^{q-1}(G/G_x)^G$. Then $\tilde \eta$
can be written $\tilde \eta=\iota_{\tilde Y} \mu_x^* \alpha $
where $\tilde Y$ is a G invariant vector field on $G/G_x$. Now
$(\mu_x)_* {\tilde Y} _{[e]} \in \kappa(Vert_xM) $ and by the
hypothesis on invariant vector fields (Theorem 2.1), there exists
a $G$ invariant vector field $Y$ on $M$ such that $Y_x = (\mu_x)_*
{\tilde Y} _{[e]}$. Thus $\tilde \eta = \mu_x^* (\iota_Y \alpha)$.

In order to calculate $d \tilde \eta$ we let $\mu$ be a $G$ basic
$(n-q)$-form which doesn't vanish at $x$ so that by
Lemma~\ref{lm53} $\alpha \wedge \mu$ is $G$ invariant $n$-form
which doesn't vanish at $x$. It is simple to check that
$\iota_\mychi [ d (\iota_Y \alpha) ]_x =0 $ if and only if
$\iota_\chi [d (\iota_Y
 \alpha) \wedge \mu ]_x=0.$

 By using the fact $\mu$ is $d$-closed and by applying equation (\ref{eq:coup}) to the invariant one-form
 $\iota_Y (\alpha \wedge \mu)= (\iota_Y \alpha) \wedge \mu$, we find
\begin{equation}
\iota_\mychi [ d (\iota_Y \alpha) \wedge \mu ]_x= \iota_\mychi [ d
( \iota _Y \alpha  \wedge \mu) ]_x  = (-1)^q \left[ d
(\iota_\mychi \iota _Y \alpha) \wedge \mu )\right]_x = 0.
\label{eq:adb}
\end{equation}
Thus $\iota_\mychi [ d (\iota_Y \alpha) ]_x =0 $.  Now computing
$$
d \tilde \eta (\tilde \mychi )_{[e]} = \left[ d (\mu_x^* \iota_Y
\alpha)(\tilde \mychi)\right]_{[e]} = \left[ \mu_x^* d (\iota _Y
\alpha) (\tilde \mychi)\right]_{[e]} = \left[ d (\iota_Y \alpha)
(\mychi) \right]_x = 0 \
$$
and using the invariance of $\tilde \eta$ we get $d \tilde \eta =
0 $. This proves ii] and therefore $H^q(\lieg, G_x) \neq 0 $.

To prove the converse, choose $x_0\in M$. Then by Theorem 4.1 the
hypothesis $A^q(\lieg ,G_x) \neq 0,$ implies there exists a
non-vanishing $G$ invariant $q$-chain $\mychi_0$ on an invariant
open neighbourhood $U$ of $x_0$. Suppose that the rank of
$\kappa(\Vertbundle)$ is $s$ and that the rank of $\kappa(TM)$ is
$r$. Let $Y_a,a=1,\ldots s$ be a local frame about $x_0$ for
$\kappa(VertM)$ consisting of invariant vector fields. Choose
invariant vector fields $Z_t,t=s+1,\ldots,r$ which together with
$Y_a$ form a local frame about $x_0$ for $\kappa(TM)$. Refine $U$
so all these objects exist on an invariant open set which we again
call $U$.

We now show that if $H^q(\lieg,G_x)\neq 0$ then $\LD_{Y_a}
\mychi_0 =0 $. First we compute
$$
(\LD_{Y_a} \alpha )(\mychi_0) =
Y_a(\alpha(\mychi_0))-\alpha(\LD_{Y_a} \mychi_0) = -
\alpha(\LD_{Y_a} \mychi_0).
$$
Expanding out the left side of this equation we get
\begin{equation}
\left[ d(\iota_{Y_a} \alpha ) + \iota_{Y_a} d\alpha \right]
(\chi_0) = - \alpha(\LD_{Y_a} \mychi_0) . \label{eq:ldv}
\end{equation}
Immediately  $ \iota_{\mychi_0} \iota_ {Y_a} d\alpha =0 $, because
$Y_a$ is vertical, while condition ii] implies by the argument
used above that$\left[ d( \iota_{Y_a} \alpha)
\right](\mychi_0)=0$, and so (\ref{eq:ldv}) along with Lemma
\ref{lm53} implies $\LD_{Y_a} \mychi_0 = 0 $.

To finish the proof of the theorem we now show there exists an
invariant $K$ which doesn't vanish at $x_0$ so that $  \mychi = K
\mychi_0$ satisfies equation (\ref{eq:ldchi}) for $Y_a,Z_t$. Using
the fact that $\LD_{Y_a} \mychi_0 = 0$, it is easy to check
$\LD_{Y_a}(K \mychi_0) = 0$. The conditions $\LD_{Z_t} (K
\mychi_0) = 0 $ leads to the differential equations for $K$
$$
\LD_{Z_t} \mychi = (Z_t(K) +K \lambda_{Z_t} )  \mychi_0= 0
$$
where $\lambda_{Z_t}$ are determined as in Lemma~\ref{lm54}.

 The functions $K, \lambda_{Z_t}$ and the vector fields
$Z_t$ are all invariant so letting $\bar Z_t = \q_* Z _t $, this
equation for $K$ can be written on $\barM$ as
\begin{equation}
\bar Z_t(\bar K) +\bar K \bar{\lambda}_{Z_t}  = 0 .
\label{eq:zbde}
\end{equation}
The integrability conditions for $K$ or $\bar K$ can be easily
verified by a computation using Lemma~\ref{lm54}. Therefore there
exists an open neighbourhood $\bar V$ of $\bar x_0$ and a
non-vanishing $\bar K$ which is a solution to (\ref{eq:zbde}).
Consequently $K \mychi_0 = (\q^* \bar K ) \mychi_0$ satisfies
(\ref{eq:ldchi}) on $\q^{-1}(\bar V)$.
\end{proof}

\section{ Examples }

\noindent {\bf Example 1.} As our first example consider the two
dimensional solvable group $G=\reals^* \times \reals$ with
coordinates $(a,b)$ acting on $\reals \times \reals^* \times
\reals$ by
$$
     (a,b)* (x,y,z) = (ax+b,ay,z).
$$
This is a free action and so $H^2(\lieg,G_x) = H^2(\lieg)$ and one
easily computes $H^2(\lieg) = 0$.

We proceed to check Theorem 5.6 for this example. The most general
$G$ invariant vertical $2$-chain $\mychi$ is of the form
$$
\mychi =K(z) y^2 \partial_x \wedge \partial_y .
$$
The invariant vector fields are
$$
R = f(z) y \partial_x +g(z) y \partial_y + h(z) \partial_z .
$$
Computing $\LD_{y \partial_y} \mychi $ we get
$$
\LD_{y \partial_y} \mychi = K(z) y^2 \partial_x \wedge \partial_y
$$
and so, consistent with Theorem 5.6 and Theorem 5.1 there is no
choice of non-zero $K(z)$ so that (\ref{eq:ldchi}) is satisfied,
and no cochain map exists.

\smallskip

\noindent {\bf Example 2.} Consider the action of the two
dimensional Abelian group $G=\reals^2$ with coordinates $(a,b)$ on
$\reals^2$ given by
$$
(a,b)* (x,y) = (x+ay+b,y) .
$$
The fact the group is Abelian implies $H^1(\lieg, G_x) \neq 0 ,$
for all $x \in \reals^2$, so a cochain map exists by Theorem 5.1.
The $G$ invariant vertical $1$-cochains are given by
\begin{equation}
\mychi = K(y) \partial _x \label{eq:ex2chi},
\end{equation}
and the invariant vector fields are
$$
R = a(y) \partial_x \ .
$$
So every cochain $\mychi $ in (\ref{eq:ex2chi}) satisfies
(\ref{eq:ldchi}). This examples demonstrates the fact that the
$q$-chain may not be unique, and by a further simple computation,
that the cochain map $\rho_\mychi$ may not be surjective.

\smallskip
As a final remark we state a theorem on the surjectivity of
$\rho_\mychi$.

\begin{theorem} Let $\mychi$ be a $G$ invariant vertical $q$-chain such that
$\rho_\mychi$ defines a cochain map. Then $\rho_\mychi$ is
surjective if and only if there exists $\alpha \in \Omega^q(M)$
such that $\alpha(\mychi) =1$ and $\alpha$ is $G$ invariant.
\end{theorem}

See \cite{and3} and \cite{felt} for other examples.

%
%
%
%

\end{document}